\documentclass[11pt]{amsart}
\usepackage[]{amsmath, amsthm, amsfonts, verbatim, amssymb}

\def\cD{{\mathcal D}}

\def\cM{{\mathcal M}}

\def\oM{\overline M}

\def\cH{\mathcal H}

\def\cT{{\mathcal T}}

\def\bC{{\mathbb C}}

\def\cT{{\mathcal T}}

\def\cM{{\mathcal M}}

\begin{document}
\newtheorem {theo}{Theorem}
\newtheorem {coro}{Corollary}
\newtheorem {lemm}{Lemma}
\newtheorem {rem}{Remark}
\newtheorem {defi}{Definition}
\newtheorem {ques}{Question}
\newtheorem {prop}{Proposition}
\def\spb{\smallpagebreak}
\def\mpb{\vskip 0.5truecm}
\def\bpb{\vskip 1truecm}
\def\wtM{\widetilde M}
\def\tM{\widetilde M}
\def\wtN{\widetilde N}
\def\tN{\widetilde N}
\def\tR{\widetilde R}
\def\tC{\widetilde C}
\def\tX{\widetilde X}
\def\tY{\widetilde Y}
\def\tE{\widetilde E}
\def\tL{\widetilde L}
\def\tQ{\widetilde Q}
\def\tS{\widetilde S}
\def\talpha{\widetilde\alpha}
\def\ti{\widetilde \iota}
\def\hM{\hat M}
\def\hq{\hat q}
\def\hR{\hat R}
\def\bs{\bigskip}
\def\ms{\medskip}
\def\ni{\noindent}
\def\td{\nabla}
\def\pd{\partial}
\def\hol{$\text{hol}\,$}
\def\Log{\mbox{Log}}
\def\bfQ{{\bf Q}}
\def\Todd{\mbox{Todd}}
\def\top{\mbox{top}}
\def\Pic{\mbox{Pic}}
\def\bP{{\bf P}}
\def\dxi{d x^i}
\def\dxj{d x^j}
\def\dyi{d y^i}
\def\dyj{d y^j}
\def\dzi{d z^I}
\def\dzj{d z^J}
\def\ozi{d{\overline z}^I}
\def\ozj{d{\overline z}^J}
\def\oz1{d{\overline z}^1}
\def\oz2{d{\overline z}^2}
\def\oz3{d{\overline z}^3}
\def\sI{\sqrt{-1}}
\def\hol{$\text{hol}\,$}
\def\ok{\overline k}
\def\ol{\overline l}
\def\oJ{\overline J}
\def\oT{\overline T}
\def\oS{\overline S}
\def\oV{\overline V}
\def\oW{\overline W}
\def\oY{\overline Y}
\def\oL{\overline L}
\def\oI{\overline I}
\def\oK{\overline K}
\def\oE{\overline E}
\def\oL{\overline L}
\def\oj{\overline j}
\def\oi{\overline i}
\def\ok{\overline k}
\def\oz{\overline z}
\def\om{\overline mu}
\def\on{\overline nu}
\def\oa{\overline \alpha}
\def\ob{\overline \beta}
\def\oGamma{\overline \Gamma}
\def\of{\overline f}
\def\og{\overline \gamma}
\def\ogamma{\overline \gamma}
\def\odelta{\overline \delta}
\def\otheta{\overline \theta}
\def\ophi{\overline \phi}
\def\opd{\overline \partial}
\def\oA{\overline A} 
\def\oB{\overline B}
\def\oC{\overline C}
\def\op{\overline D}
\def\ov{\overline \varphi}
\def\oIq1{\oI_1\cdots\oI_{q-1}}
\def\oIq2{\oI_1\cdots\oI_{q-2}}
\def\op{\overline \partial}
\def\ua{{\underline {a}}}
\def\us{{\underline {\sigma}}}
\def\Chow{{\mbox{Chow}}}
\def\vol{{\mbox{vol}}}
\def\dim{{\mbox{dim}}}
\def\rank{{\mbox{rank}}}
\def\diag{{\mbox{diag}}}
\def\tor{\mbox{tor}}
\def\bp{{\bf p}}
\def\bk{{\bf k}}
\def\a{{\alpha}}
\def\tchi{\widetilde{\chi}}
\def\ta{\widetilde{\alpha}}
\def\ovarphi{\overline \varphi}
\def\ocH{\overline{\cH}}
\def\tV{\widetilde{V}}
\def\tf{\widetilde{f}}
\def\th{\widetilde{h}}

\ni
\title[Tower of coverings of quasi-projective varieties]
{Tower of coverings of quasi-projective varieties}
\author[Sai-Kee Yeung]
{Sai-Kee Yeung}

\address[]{Mathematics Department, Purdue University, West Lafayette, IN  47907
USA} \email{yeung@math.purdue.edu}

\thanks{ 
The author was partially supported by a grant from the National Science Foundation}

\ni{\it }

\maketitle
\ni{\bf Abstract} {\it The main goal of this article is to relate asymptotic geometric properties on a tower of coverings
of a non-compact K\"ahler manifold of finite volume with reasonable geometric assumptions
to its universal covering.
 Applicable examples 
 include moduli spaces of hyperbolic punctured Riemann surfaces and Hermitian locally symmetric spaces of finite
volume.} 
\bigskip
\noindent

\bigskip
\begin{center}
{\bf 1. Introduction}
\end{center}

\ms
\ni {\bf 1.1}  Let $M$ be a complex manifold.   By {\it a tower of coverings} of $M,$ we mean a sequence of  finite coverings $M_{i+1}\rightarrow M_i$
with $M_0=M,$ such that $\pi_1(M_{i+1})$ is a normal subgroup of $\pi_1(M_1)$ of finite index and $\cap_{i=0}^\infty\pi_1(M_i)=\{1\}.$ 
An interesting problem is to relate the geometric properties of $M_i$ to $\tM,$ 
the universal covering of $M.$ 
The case $M$ is compact has been an object of study for a long time, cf. [DW], [Do], [L], [K1], [O1], [R], [SX], [X], [To], [Y1], [Y3]
and many more.

The study for a general non-compact $M$ of finite volume
with respect to some complete K\"ahler metric has been limited.      Since many interesting geometric and arithmetic objects arise as non-compact complex manifolds, it is natural and meaningful to ask if similar results hold for non-compact manifolds under mild restrictions.   In fact, the paper grows out of curiosity to understand
the corresponding asymptotic behavior for general moduli spaces of hyperbolic punctured Riemann surfaces. The other motivation is to understand
if asymptotic properties satisfied by compact Hermitian locally symmetric spaces of non-compact type as studied in [D], [To], [Y1] and [Y3]  
are also satisfied by non-compact ones.  For non-compact Hermitian locally symmetric spaces of finite volume, results about asymptotic growth of Betti numbers similar to the results of [DW]
have been obtained by [BM], [De], [K2], [S], but results related to [Do], [To], [Y1], or [Y3]  appear to be open up to now.  Our main goal 
in this paper is to present a formulation which is applicable to a general non-compact manifold with
reasonable restrictions which in particular include the two classes of manifolds mentioned above.

The usual difficulty to discuss asymptotic growth of geometric quantities such as cohomology
 on a tower of non-compact manifolds
is that in general a big proportion of the quantity may escape 
to the infinity as one takes an appropriate limit.   We show that the difficulty can be overcome under reasonable conditions on the manifolds.
Specifically, we assume that the manifolds involved are 
 geometrically finite or
quasi-projective, to be explained in {\bf 2.1}.  These conditions are natural and are satisfied by both
moduli spaces of hyperbolic Riemann surfaces and Hermitian locally symmetric spaces of finite volume.

As a result, we are able to generalize estimates from towers of
compact manifolds to similar towers of non-compact manifolds, such as relation between
growth
of Betti numbers with respect to volume and the $L^2$ Betti numbers (von Neumann dimensions),
 convergence of Bergman kernels, and equidistributions of
pluricanonical sections.  We also verify that the canonical sections of a sufficiently 
large cover in a tower of quasi-projective varieties give rise to an immersion of the manifold into some projective space.

It is a pleasure for the author to thank Gopal Prasad for helpful comments on the exposition of
the paper.  The author is also grateful to the referee for pointing out numerous misprints in an early version of this paper and providing constructive 
suggestions.

\bigskip
\begin{center}
{\bf 2. Statement of results}
\end{center}

\ms
\ni {\bf 2.1} Suppose $M$ is a complex manifold 
equipped with a K\"ahler metric $g$ of finite volume.
We denote the K\"ahler form of $g$ by $\omega.$  Let $\tM$ be the universal covering of
$M.$  

\ms
\ni{\bf Definition} We say that $(M,g)$ is {\it geometrically finite}
if \\
(i) the volume of $M$ with respect to $g$ is finite,\\
(ii) the Riemannian sectional curvature of $g$ is uniformly bounded from above, and\\
(iii) the injectivity radius of $\tM$ is uniformly bounded from below on $\tM$.

\ms 
As usual, we say that $M$ is {\it quasi-projective} if it can be written as $M=\oM-D,$ where $\oM$ is a projective algebraic
manifold, and $D$ is a divisor on $\oM.$  In case $M$ is quasi-projective, using resolution of singularities 
we can choose $\oM$ 
in such a way that  $D=\oM-M$ is a divisor with normal crossings.  Denote by $\Delta$ the
unit disk in $\bC$ and $\Delta^*$ the punctured unit disk.  We may cover a neighborhood of $D$ in
$\oM$ by a finite number of open sets of the form $U=\Delta^{n-k}\times (\Delta^*)^{k},$ where
$n\geqslant k\geqslant 1.$  Clearly $\oM$ is covered by a finite number of $U$ as above if we allow
$n\geqslant k\geqslant 0.$

\ms
\ni {\bf 2.2}  Suppose that there
exists a tower of coverings of $M$ as mentioned in {\bf 1.1}.  We
 let $\Gamma=\pi_1(M)$ be the fundamental group of
$M.$
Let $\Gamma_1=\Gamma$ and $\Gamma_1<\dots<\{1\}$ be a tower of normal subgroups of $\Gamma$ corresponding
to an infinite sequence of normal coverings with finite index of $M,$  so that $\cap_{i=0}^\infty\Gamma_i=\{1\}.$  
In other words, the fundamental group is residually finite.  Let $D_i$ be a fundamental domain of $\Gamma_i.$  Since we are taking a tower of
  normal coverings, we may assume that the fundamental domains $D_i$ of $\Gamma_i$ are nested in the sense that $D_i\subset D_{i+1}$.  As
  $\cap_i \Gamma_i=1,$ $\tM=\cup_i D_i.$ 
 If $M_1$ is quasi-projective, each $M_j$ for $j>1$ is quasi-projective as well. 
 
 \ms
\ni {\bf 2.3}   Denote by
$H_{(2)}^{p}(N)$ the space of  $L^2$ harmonic $p$-forms on a complete Riemannian manifold $N,$  with its dimension 
denoted by $b^p_{(2)}(N).$  If moreover, $N$ is a complete K\"ahler manifold, we let 
$H_{(2)}^{p,0}(N)$ be the space of $L^2$ holomorphic $p$-forms on a manifold $N,$ which by K\"ahler
identity, is isomorphic to $H_{(2)}^{0,p}(N).$  We denote by
$h_{(2)}^{p,0}(N)$ its dimension.

Let $\cD$ be a fundamental
domain of $\pi_1(N)$ in $\tN.$  The {\it von Neumann 
dimension} of $H_{(2)}^{p,0}(\tN)$ with respect to $N$ is defined as 
$\int_{\cD} B^{p,0}_{\tN}(x),$ where $B^{p,0}_{\tN}(x)=B^{p,0}_{\tN}(x,x)$ is the trace of the Bergman kernel as defined in
{\bf 3.1}, and is denoted by $h_{v,(2)}^{p,0}(\tN).$   
We may also consider them as the dimension of the corresponding Dolbeault
cohomology from Leray Isomorphism.  Similarly, we define the von Neumann dimension of the space
of $L^2$ $d$-harmonic forms $p$ on $\tN,$ and denote the dimension by $b_{v,(2)}^{p}(\tN).$

\ms
\ni {\bf 2.4}  To relate geometric properties of a tower of quasi-projective manifolds
to its universal covering, we begin with some qualitative asymptotic statements.

\ms
\ni{\bf Definition} We say that a complete manifold $N$ of complex dimension $n$ satisfies
{\it cohomology condition C} if $H_{(2)}^{p,0}(\tM)=0$ for $p<n$ and $H_{(2)}^{n,0}(\tM)\neq0.$

\ms

 For a tower of compact K\"ahler manifolds, the following result is Theorem 1.3 of [Do], see also [R] for the case of Riemann surfaces. 
 It also follows immediately from
 Theorem 1.1 of [Y1], since part (b) is an immediate consequence of part (a) as illustrated by the last sentence in {\bf 3.4}.   The interest in this article is the non-compact version.  
  
 \begin{theo}
Let $M$ be a complex manifold of complex dimension $n$ equipped with a K\"ahler metric which is geometrically
finite.  Assume that $M$ supports
a tower of normal coverings $M_i$ of $M.$   Assume that $\tM$ satisfies the cohomological condition C.  Then for each $0\leqslant p\leqslant n$,\\
(a). 
\begin{equation}
\lim_{i\rightarrow\infty}\frac{h_{(2)}^{p,0}(M_i)}{[\pi_1(M):\pi_1(M_i)]}=h_{v,(2)}^{p,0}(\tM).
\end{equation}
(b). For each point $x\in D\subset D_i\subset \tM,$ the Bergman kernels satisfy
\begin{equation}
\lim_{i\rightarrow\infty} B^{p,0}_{M_i}(x)=B^{p,0}_{\tM}(x).
\end{equation}
\end{theo}

\ms
\ni{\bf Remarks}

\ms
\ni{\bf (a).} As mentioned earlier, the result for $M$ compact is already known in various cases, (cf. for instance [Y1]).  From
this point onward, we will focus on the non-compact cases.  

\ms
\ni{\bf (b).} The existence of a tower of coverings on a manifold is
not always available.   On the other hand,
there are lots of natural examples supporting such towers and moreover satisfying other conditions stated in the theorem, including
Hermitian locally symmetric spaces, moduli space of Riemann surfaces with punctures which are hyperbolic, and manifolds with
non-positive sectional curvature.  We refer the readers to Section {\bf 4} for the details.  Furthermore we only focus
on normal coverings of $M.$  We refer the readers to [O2] for discussions about necessity of such conditions.

\ms
\ni{\bf (c).}  As will be clear from the proof, the same argument also implies that for a tower of (non-compact) manifolds $M_i$ for which
the universal covering $\tM$ satisfies
$b^p_{v,(2)}(\tM)=0$ for $p\neq n,$ we conclude that 
$$\lim_{i\rightarrow\infty}\frac{b^p_{(2)}(M_i)}{[\pi_1(M):\pi_1(M_i)]}=b_{v,(2)}^{p}(\tM)$$
for all $p.$

\ms
\ni{\bf (d)}. In the case of compact Hermitian locally symmetric spaces of non-compact type,  the paper of Kazhdan in [K1]
gives the first result observing that 
$\lim_{i\rightarrow\infty}\frac{h_{(2)}^{n,0}(M_i)}{[\pi_1(M):\pi_1(M_i)]}>0,$  which
was later proved also for non-compact Hermitian locally symmetric spaces by Kazhdan [K2], see also [S], page 149, Corollary 1.  Theorem 1(a) and (b)
can be considered to be a more precise version of the above results and is applicable to examples such as moduli spaces of curves
with punctures, cf. {\bf 5.2}.

\ms
\ni {\bf 2.5} From the point of view of automorphic forms or cusps forms, the following result may be interesting.
\begin{theo}
Let $M=\oM-D$ be a quasi-projective variety equipped with a K\"ahler metric which is geometrically finite.
Assume that $M$ supports 
a tower of normal coverings $M_i$ of $M.$  Let $K$ be the canonical line bundle on $M.$  Let $L$ is a positively
curved Hermitian line bundle on $M$.   Then $\lim_{i\rightarrow\infty}\frac{h_{(2)}^{0}(M_i,K+L)}{[\pi_1(M):\pi_1(M_i)]}=h_{v,(2)}^{0}(\tM,K+L).$\\
\end{theo}

\ms\ni{\bf Remarks}

\ms
\ni{\bf (a).}  Classical automorphic forms on Hermitian symmetric spaces of non-compact type correspond
to $L=\ell K$ or $K+L=(\ell + 1) K,$ where $\ell$ is a positive rational number for which $\ell K$ is a line
bundle on $M_i$ for $i$ sufficiently large.

\ms
\ni{\bf (b).}  The same formulation is applicable to cusp forms, which can be considered as the space of 
pluri-logarithmic canonical forms vanishing at the compactifying divisor. 

\ms
\ni{\bf (c).}  The proof of Theorem 2 for cusp forms for non-compact Hermitian locally symmetric spaces of finite
volume has been given in various settings in the work of [K2], [BM], [De] and [S].

\ms
\ni {\bf 2.6} Recall that the Bergman metric of a complex manifold can be defined as
$\sqrt{-1}\pd\op\log B_M,$  where $B_M=B_M^{n,0}.$    For a general complex manifold,
the $(1,1)$ form gives rise to a pseudo-metric which may not be positive definite.

The following is a result to recover the Killing metric on a Hermitian symmetric space.

\begin{theo}
Let $M$ be a quasi-projective variety supporting a tower of coverings as studied in Theorem 2.
Then $\sqrt{-1}\pd\op\log B_{M_i}$ converges on compacta to
the Bergman metric on $M.$ 
\end{theo}

In case $M$ is a Hermitian locally symmetric space, the Bergman metric
is just the invariant Killing metric up to a normalizing constant. 
For the special case that $M$ is a compact hyperbolic Riemann surface, the above theorem was a theorem of 
Rhodes [R].  Donnelly [Do] generalized the result to a compact Hermitian locally symmetric space
of non-compact type.  The theorem here covers non-compact Hermitian locally symmetric spaces of finite
volume as well as moduli spaces of curves with punctures.  

\ms
\ni{\bf 2.7}  We would consider two applications of the earlier results to Hermitian locally symmetric spaces
and moduli spaces of curves.  In both cases, the universal covering is biholomorphic to a bounded domain in
some $\bC^n.$
First of all, Theorem 3 leads immediately to
equidistribution of a generic $L^2$-section of the canonical line bundle.
The result for a compact tower has been achieved by To in [To],  see also [SZ] for formulation in related
directions.

\begin{theo}  Suppose that we are given a tower of coverings of Hermitian locally symmetric spaces
or moduli space of punctured Riemann surfaces as discussed earlier.   Denote by $Z_s$ the current
of integration associated to the zero
divisor of $s\in H^0_{(2)}(M_i, K_{M_i}).$  We may regard $Z_s$ as a random variable as $s$ varies over the set of holomorphic
sections of $H^0_{(2)}(M_i,K_{M_i})$ with $L^2$ norm $1.$  We refer the readers to {\bf 4.3} for more details on the settings.
The expected values of $Z_s$
satisfies
$$\lim_{i\rightarrow\infty}{E_i}(Z_s)=\frac1{2\pi}\sqrt{-1}\pd\op\log B_{\tM}.$$
\end{theo}

\ms
The second application is the following.  Similar to the results in [Y3] for cocompact lattices of Hermitian locally symmetric spaces, we
have the following consequence.

\begin{theo}
Let $M_i$ be a tower of coverings as studied in Theorem 4.
There exists $i_o>0$ such that for $i\geqslant i_o,$ global $L_2$ holomorphic sections in $\Gamma(M_i,K_{M_i})$ 
give rise to a holomorphic immersion of 
$M_i$ into some projective space.
\end{theo}

Note that the interest here is that sections of $K_{M_i}$ instead of multiples of $K_{M_i}$ give immersion of 
$M_i.$

\bigskip
\begin{center}
{\bf 3. Bergman kernels and asymptotic results.}
\end{center}

\ms
\ni {\bf 3.1} In this section, $M$ is a complex manifold of complex dimension $n$ as studied in Section 2. 
Let us recall some standard terminologies.
Let $\varphi\in H_{(2)}^{p,0}(M).$   The $L^2$-norm of $\varphi$ is defined by
\begin{eqnarray}
\Vert\varphi\Vert^2&=&\int_M \varphi\wedge*\varphi,
\end{eqnarray}
which can also be expressed as
\begin{equation}
\int_M|\varphi|^2=\int_M|\varphi|_g^2\omega^n,
\end{equation}
where $|\cdot|_g$ is the norm
with respect to the K\"ahler metric $g$ associated to $\omega,$ and $\omega^n$ is the volume form of the metric $g$ on $M.$  

Let $\{f_i\}$ be an orthonormal basis of $H_{(2)}^{p,0}(M).$  The Bergman
kernel is defined to be 
$$B_M^{p,0}(x,y)=\sum_i f_i(x)\wedge *f_i(y),$$
where $*$ is the Hodge operator.  As such we are regarding $B_M^{p,0}$ as a section of
$p_1^*\Omega_M^{p,0}\otimes p_2^*\Omega_{M}^{n-p,n},$ where $p_i$ is the projection of
$M\times M$ into the $i$-th factor.

We are mainly interested in the trace of the Kernel, $B_M^{p,0}(x):=B^{p,0}_M(x,x).$
From definition, 
\begin{equation} \label{ex}
h^{p,0}_{(2)}=\int_M B_M^{p,0}(x,x).
\end{equation}

As the Bergman kernel is independent of basis, for each fixed point $x\in M,$
$$B_M^{p,0}(x,x)=(\sup_{f\in H^{p,0}_{(2)}(M),\Vert f\Vert=1}|f(x)|_g^2)\omega^n=\sup_{f\in H^{p,0}_{(2)}(M),\Vert f\Vert=1}|f(x)|^2,$$ 
Here by abuse of notations, we denote $|f(x)|_g^2\omega^n$ by $|f(x)|^2.$  Note that the space of $(n,n)$ forms are point wise
one dimensional.  In the same manner, for two $(p,0)$ forms $f_1$ and $f_2,$ we say that 
$|f_1|^2=f_1\wedge*f_1\leqslant |f_2|^2=f_2\wedge*f_2$ if $|f_1|^2_g\leqslant |f_2|_g^2.$

For $p=n,$ we may also write
$$B_M^{n,0}(x,x)=\sup_{f\in H^0_{(2)}(M,K_M),\Vert f\Vert=1}|f(x)|^2.$$

\ms
\ni {\bf 3.2}   In this section, we do not need to assume that the manifold $M$ involved is 
quasi-projective, but assume that the K\"ahler metric involved is geometrically finite
and the universal covering satisfies cohomology condition C.  We begin with the following
observation.

\begin{lemm}\label{lem3}
Let $x\in M.$  Identify $x$ with a point (still denoted by) $x\in D\subset D_i$ on the universal
covering $\tM$ for all $i.$  Then for $0\leqslant p\leqslant n,$ where $n=\dim_{\bC}M,$ we have
$$\limsup_{i\rightarrow\infty} B^{p,0}_{M_i}(x)\leqslant B^{p,0}_{\tM}(x).$$
\end{lemm}

\ni{\bf Proof}   We may assume that $B_{M_i}^{(p,0)}(x)$ is realized by $|f_{i,x}(x)|^2$ for some 
$f_{i,x}\in H^{p,0}_{(2)}(M_i)$ with $\Vert f_{i,x}\Vert_{L^2(M_i)}=1.$  Suppose $\limsup_{i\rightarrow\infty}|f_{i,x}(x)|= A.$
Consider the sequence of forms $f_{i,x}.$  We are going to show that $f_{i,x}$ converges on compacta to $f_{\infty, x}$ with $L^2$ norm bounded from above by $1.$

Let $V$ be any relatively 
compact set of $\tM.$  Note that we may assume that $D_{i+1}\subset D_{i}$.  Since $\cup_i D_i=\tM$ from definition, 
we conclude that $V\subset D_i$ for all $i$ sufficiently large.

Now we claim that by taking a subsequence if necessary, the sequence $f_{i,x}\vert_V$ is equicontinuous.

Let $\chi_{D_i}$ be the 
characteristic function on $D_i.$  By considering $\chi_{D_i}\tf_{i,J,x}$ we may regard
$\tf_{i,x}$ as a function on $D_i.$  By taking a subsequence if necessary, we know that as elements in the Hilbert space $H^{p,0}_{(2)}(D_i),$
 $\tf_{i,x}$ form a Cauchy sequence as
$i\rightarrow\infty$ for fixed $x.$  In particular, given any $\epsilon>0,$ there exists $N>0$ such that
$\Vert \tf_{k,x}-\tf_{j,x}\Vert_{D} \leqslant \epsilon$
if $k\geqslant j\geqslant N.$   The same conclusion holds when $D$ is replaced by $D_i.$
As $V\subset D_i$ for $i$ sufficiently large, 
we conclude that the $L^2$ norm of $(\tf_{k,x}-\tf_{j,x})\vert_V$ is bounded by $\epsilon<1.$ 

Since $\tM$ has bounded geometry, the norm square of a harmonic $(p,0)$ form $\varphi$ satisfies subelliptic differential inequality
of the form
$$\Delta|\varphi|_g^2+k|\varphi|_g^2\geqslant 0,$$
(cf. page 204 of [Y1]).  Then standard regularity theory implies that the point wise norm of a form $\varphi$ is bounded by the $L^2$ norm of $\varphi$ 
as given in page 203-205 of
[Y1], from which we conclude that the pointwise norm of $\tf_{i,J,x}$ is equicontinuous on
$V.$  Hence the claim follows.

An alternate simple argument of the above paragraph is as follows.  Let $U_x$ be a small complex coordinate neighborhood
of $x$ on $\tM.$  The earlier argument implies that  the $L^2$ norm of $(\tf_{k,x}-\tf_{j,x})\vert_{U_x}$ is bounded by $\epsilon<1.$ 
On $U_x,$ 
a holomorphic $(p,0)$-form can be written
as $\sum_{J}\tf_{i,J,x}dz^{j_1}\wedge\cdots\wedge dz^{j_p}$ in terms of local coordinates, where the sum is over all
$p$-tuples $J=(j_1,\cdots, j_p)$ with $j_1<\cdots <j_p.$
Hence convergence properties of  $f_{i,x}$ on compacta is the same as the convergence of
the local holomorphic function $\tf_{i,J,x}.$ 
With the knowledge of $L^2$ bounds, the Maximum
Principle (or Cauchy Estimate) implies the point wise estimate of $|\tf_{k,J,x}-\tf_{j,J,x}|$ as well.  Since $V$ as 
a relatively compact set can be covered by a finite number of such neighborhoods, 
the claim follows.

From the claim, we apply Ascoli-Arzela Theorem to conclude that given any sufficiently small $\epsilon>0,$ there
exists a subsequence of
$f_{i_x}$ converges on compacta to a holomorphic form $f_x\in H^{p,0}_{(2)}(\tM),$
with $|f_x(x)|\geqslant A-\epsilon$ and 
$$\Vert f_x\Vert_V\leqslant \limsup_{i\rightarrow\infty}\Vert f_{i,x}\Vert_V\leqslant 1.$$
Since $V$ is an arbitrary compact subset of $\tM,$ by considering an nested exhaustive sequence of such
$V$, a standard argument involving normal family of functions 
 concludes the proof of Lemma \ref{lem3}.

\qed

\ms
\ni {\bf 3.3}  Let us recall that on a complete K\"ahler manifold $M$ of complex
dimension $n$ with finite volume,  the $L^2$-arithmetic genus  $\chi_{(2)}(M)$ and
$L^2$ Euler-Poincar\'e characteristic $e_{(2)}(M)$ are defined by
\begin{eqnarray*}
\chi_{(2)}(M)&=&\sum_{p=0}^n(-1)^ph_{(2)}^{p,0}(M),\\
e_{(2)}(M)&=&\sum_{j=0}^n(-1)^pb_{(2)}^j(M),\\
\end{eqnarray*}
when the expressions involved are finite.
Similarly, we define the corresponding von-Neuman dimension of
$\tM$ with respect to $\pi_1(M)$ by
\begin{eqnarray*}
\chi_{v,(2)}(\tM)&=&\sum_{p=0}^n(-1)^ph_{v,(2)}^{p,0}(\tM),\\
e_{v,(2)}(\tM)&=&\sum_{j=0}^n(-1)^pb_{v,(2)}^j(\tM),\\
\end{eqnarray*}

\begin{lemm}  \label{lem4}
Assume that $(M,g)$ is geometrically finite.  Then \\
(i). The arithmetic genus satisfies $\chi_{v,(2)}(\tM)=\chi_{(2)}(M)$\\
(ii).  The Euler-Poincar\'e number satisfies $e_{v,(2)}(\tM)=e_{(2)}(M)=e(M)$
\end{lemm}
\ni{\bf Proof}  in case $M$ are compact, this is just Atiyah's Covering Index Theorem [A].
For $M$ non-compact of finite volume 
 and is geometrically finite, the results are still valid as observed by Cheeger and Gromov
in [CG].  The idea is to make use of a good exhaustion of the manifold, on which the curvature and the second fundamental forms of the 
boundary of the exhaustion can be estimated.   Geometric
finiteness properties of the K\"ahler metric is used to construct a good exhaustion.   Then the usual proof of the
Atiyah's Covering Index Theorem in terms of the traces of the heat kernels of differential forms can be adopted to this case, as 
given in Section 6 of [CG].

\qed

\ms
\ni {\bf 3.4} {\bf Proof of Theorem 1}
Note that from definition, 
$$\int_{M_i}B_{M_i}^{p,0}(x)=h^{p,0}_{(2)}(M_i).$$
Since the Bergman kernel is invariant under biholomorphism and the coverings involved
are normal coverings, the left hand side can be expressed as
$$[\Gamma,\Gamma_i]\int_M B_{M_i}^{p,0}(x)=[\Gamma,\Gamma_i]\int_D B_{M_i}^{p,0}(x).$$
From Lemma \ref{lem3}, we conclude that 
\begin{eqnarray*}
h^{p,0}_{v,(2)}(\tM)&=&\int_D B_{\tM}^{p,0}(x) \\
&\geqslant&\limsup_{i\rightarrow\infty} \int_D  B_{M_i}^{p,0}(x)\\
&=&\limsup_{i\rightarrow\infty} \frac{h^{p,0}_{(2)}(M_i)}{[\Gamma,\Gamma_i]}
\end{eqnarray*}

Now for $p<n,$ we know from assumption that $B_{\tM}^{p,0}=0$ and hence
$h^{p,0}_{v,(2)}(\tM)=0.$  It follows that $\limsup_{i\rightarrow\infty} \frac{h^{p,0}_{(2)}(M_i)}{[\Gamma,\Gamma_i]}=0.$
Hence automatically
\begin{equation}\label{small}
\lim_{i\rightarrow\infty} \frac{h^{p,0}_{(2)}(M_i)}{[\Gamma,\Gamma_i]}=h^{p,0}_{v,(2)}(\tM)
\ \ \ \mbox{for $p<n.$}
\end{equation}

On the other hand, from Lemma \ref{lem4}, for each $i,$
$\chi_{v,(2)}(\tM)=\frac{\chi_{(2)}(M_i)}{[\Gamma,\Gamma_i]},$ which implies that

$$\sum_{p=0}^n(-1)^p h_{v,(2)}^{p,0}(\tM)=\sum_{p=0}^n(-1)^p \frac{h_{(2)}^{p,0}(M_i)}{[\Gamma,\Gamma_i]}.$$

From equation \ref{small}, after taking limit as $i\rightarrow\infty,$ this implies that
$\lim_{i\rightarrow\infty} \frac{h^{n,0}_{(2)}(M_i)}{[\Gamma,\Gamma_i]}=h^{n,0}_{v,(2)}(\tM).$
This concludes (i) of Theorem 1.  (ii) follows from (i) and the interpretation in (\ref{ex}).

\qed

\ms
\ni {\bf 3.5} {\bf Proof of Theorem 2}   The idea of proof is similar to the one in Theorem 1.
The K\"ahler metric $g$ on $M$ induces a Hermitian metric $g_K$ on the canonical line bundle
$K$ of $M$.  Note that
$g_K$ is just the reciprocal of the determinant of the K\"ahler metric in terms of local
coordinates.  Let $h$ be a positively curved metric on $L.$  One defines the Bergman
kernel of $K+L$ on $M$ as $B_{M, K+L}$ at $x\in M$ as $\sum_i|f_i|^2_{g,h}(x),$ where $\{f_i\}$
is an orthonormal basis of $K+L$ with respect to metrics $g_K\cdot h$ on $K+L$ and the volume
form of $g$ on $M.$   Similarly,  
$|\cdot|_{g,h}$ denotes the pointwise norm of the section with respect to $g_K$ and $h.$

Similar to the case of differential forms, we denote the dimension of the space of $L^2$ holomorphic $K+L$ valued forms on $M$ by
$h^{p,0}_{(2)}(M, K+L)$ and the von Neumann dimension of $K+L$ with respect to $M$ by
$h^{p,0}_{v,(2)}(\tM, K+L)=\int_D B_{\tM, K+L}^{p,0}(x),$ where $D$ is a fundamental domain
of $\pi(M)$ in $\tM.$

Since $(L,h)$ is a positively curved Hermitian line bundle, Kodaira Vanishing Theorem
implies the vanishing of $h^{p,0}_{(2)}(M, K+L)$ and $h^{p,0}_{(2)}(\tM, K+L).$ 
The latter implies the vanishing of $B_{\tM, K+L}^{p,0}$ and hence the
vanishing of $h^{p,0}_{v,(2)}(\tM, K+L).$  The rest of the proof is then the same as Theorem 1.
Again, the problem of non-compactness is overcome since only $L^2$ sections are concerned, and the good exhaustion of Cheeger-Gromov [CG] can be applied to complete the argument
using Atiyah's Covering Index Theorem.

\qed

\ms
\ni {\bf 3.6} {\bf Proof of Theorem 3}  From Theorem 1(ii), we conclude that
$B_{M_i}^{p,0}$ converges pointwise on compacta to $B_{\tM}^{p,0}.$  Note that each
$B_{M_i}^{p,0}(x,y)$ as well as $B_{\tM}^{p,0}$ expressed
in terms of local coordinates in a coordinate neighborhood is analytic as a function on $M_i\times\overline{M_i},$ where
$\overline{M_i}$ is the complex manifold whose underlying differentiable structure
is the same as $M_i$ but the complex structure is the complex conjugate.  The argument
of Theorem 1 clearly also shows that $B_{M_i}^{p,0}$ converges in $C^k$ to $B_{\tM}^{p,0}$ 
for all $k$.   Theorem 3 follows from $C^2$ convergence.
\qed

\ms
\ni {\bf 3.7} We would like to give a few remarks.

\ms
\ni{\bf Remarks}

\ms
\ni{\bf (a)}   By considering Taylor series expansion, and note that in terms of local coordinates, the coefficients of $B_{M_i}^{p,0}(x,y)$
in terms of standard basis for the differential forms
is holomorphic in $x$ but antiholomorphic in $y,$ the convergence along the diagonal
of $M\times \oM$ implies convergence everywhere on $M\times\oM.$  Hence one
actually has analytic convergence on $M\times\oM$ and hence along the diagonal as well.  We refer the readers to [Y3] for
details of the arguments.

\ms
\ni{\bf (b)} As mentioned in the introduction, Kazhdan proved in [K1] for compact 
Hermitian locally symmetric spaces that $h_{v,(2)}^{n,0}(M_i)>0,$  or that the Bergman kernel on 
$\tM$ is positive.
This was utilized to prove
the following important result.  Let us call a quotient of Hermitian symmetric space 
by a co-compact arithmetic lattice, an arithmetic variety, which is known to be defined over
some number field, cf. [K1].  As a
variety defined over a number field, it is interesting to study the conjugate of an arithmetic variety by an element in the absolute Galois group.  Kazhdan
proved in [K1] that such a conjugate is also an arithmetic variety, in other words, another arithmetic
quotient of a Hermitian symmetric space.
For non-cocompact lattices, proofs have been given by [K2] and
Nori and Raghunathan [NR].  A completely different geometric proof for all cases
has been given in Mok-Yeung [MY].  
The readers may also consult [Mi] for more exposition and remarks on this result of Kazhdan.

\ms
\ni {\bf 3.8} Let us now explain that the results of Theorem 1 are also applicable to
harmonic forms and the usual Betti numbers.  In the first place, Lemma \ref{lem3} is applicable to
harmonic $(p,q)$ forms.  
The reason is that a harmonic form satisfies an elliptic equation
which becomes uniformly elliptic on a relatively compact set $V.$ 
Lemma \ref{lem4} is also applicable for such harmonic forms, by considering $\chi^{(q)}=\sum_{p=0}^n(-1)^ph_{(2)}^{p,q}(M).$  The rest of the argument is the same as in the proof of
Theorem 2.

\bs

\begin{center}
{\bf 4. Hermitian locally symmetric spaces and moduli space of hyperbolic punctured Riemann surfaces.}
\end{center}

\ms
\ni {\bf 4.1} In this section, we explain briefly the reason that non-compact
Hermitian locally symmetric
spaces of non-compact type 
and moduli space $\cM_{g,n}$ of  hyperbolic Riemann surfaces of genus $g$ with $n$ punctures satisfy the hypothesis required for
our theorems.

A Hermitian locally symmetric space can be written as $M=\Gamma\backslash G/K,$ where $G$ is a semi-simple
Lie group, $K$ is a maximal compact subgroup and $\Gamma$ is a lattice so that $M$ has finite
volume with respect to the invariant metric, the Bergman metric.  For simplicity, we may just consider
a torsion-free lattice.  Otherwise we have to consider some \'etale coverings to resolve the singularities
which are quotient singularities.  The Bergman metric on such manifolds has non-positive 
Riemannian sectional curvature and the volume is finite.  
It is well-known that $\Gamma$ is
residually finite.  Explicit examples are given by arithmetic lattices.  A tower of coverings is
then obtained by considering a tower of congruence subgroups of $\Gamma.$

Since the universal covering of such a manifold is a bounded symmetric domain, it is
clear that the geometry is finite.  Moreover, it is also well known that the Bergman metric
satisfies Condition C, (cf. [DW], [BW] or [DF]).

\ms
\ni {\bf 4.2} For the moduli space of Riemann surfaces $\cM_{g,n},$ it is known that the mapping
class group $\Gamma_{g,n}$ of $\cM_{g,n}$ is residually
finite according to a result of Grossman (cf. [G]).  Hence there exists a tower of normal subgroups $\Gamma_i$ with
$\Gamma_1=\Gamma_{g,n}$ and $\cap_{i=1}^\infty \Gamma_i=\{1\}.$  The universal covering of $\cM_{g,n}$ is the Teichm\"uller
space $\cT_{g,n}.$  The tower of normal coverings are denoted by $\cT_{g,n}/\Gamma_i.$

In general, $\cM_{g,n}$ contains quotient singularities corresponding to the fixed points of the mapping class
group.  The singularities can be resolved by considering level structure.  We may consider such a finite
normal covering to begin our study.  We refer the readers to [HM] for backgrounds on moduli space of curves.

It is also known that $\cM_{g,n}$ supports a
Bergman metric for which the moduli space, or the Teichm\"uller space,
satisfies Condition C and has finite geometry.
See for example [Y5] or the earlier work in [Mc].

\ms
\ni {\bf 4.3}  We now elaborate on the setting and the proof of Theorem 4.   The setting is similar to the one 
given by To in [To], see also [SZ].   It is known that the space of $L^2$ sections on $H^0_{(2)}(M_i,K_{M_i})$ is of finite dimension.
This follows for example from a well-known argument of Siegel (cf. [Mo1]).
Let $SH^0_{(2)}(M_i,K_{M_i})$ be the set of holomorphic sections of the canonical line bundle with $L^2$ norm $1$ on $M_i,$
equipped with the standard Haar measure $\mu_i.$  Denote by 
$\cD^{1,1}(M_i)$ the space of $(1,1)$-currents on $M_i.$  The divisor of any $s\in H^0_{(2)}(M_i,K_{M_i})$
defines a current $Z_s\in \cD^{1,1}(M_i)$.  Then as $s$ varies over the probability space
$(SH^0_{(2)}(M_i,K_{M_i}),\mu_i),$ $Z_s$ can be regarded as a $\cD^{1,1}(M_s)$ valued 
random variable.  The expectation value $E_i(Z_s)\in \cD^{1,1}(M_i)$ is defined by
$$(E_i(Z_s),\eta):=\int_{s\in SH^0_{(2)}(M_i,K_{M_i})}\big(\int_{Z_s}\eta)d\mu_i(s),$$
for any test smooth $(1,1)$ form $\eta.$  Since we are considering a normal tower of coverings,
the expected values are invariant under deck-transformations and we may just regard it as
living on $M$ or its fundamental domain.

\newpage
\ni{\bf Proof of Theorem 4} Once we have Theorem 3, the argument of To in [To] can immediately
be applied to conclude the proof of Corollary 1.  The argument is related to the arguments of
Shiffman and Zeltditch in [SZ].

\qed

\bs
\ni {\bf 4.4} {\bf Proof of Theorem 5}   Equipped with Theorem 3, the scheme of proof is similar to the cocompact case as in
[Y3], [Y4].  However, we need to pay attention to the fact that we are considering non-compact manifolds in which the injectivity radius at a point $x$ approaches zero as $x$ approaches
boundary of the manifold.    We begin with the following observation.

\begin{lemm}\label{lem1}
(i). A $L^2$ holomorphic $n$-form on a quasi-projective $M$ can be extended as a holomorphic $n$-form to $\oM.$
\end{lemm}

\ni{\bf Proof}  The $L^2$ norm of a holomorphic $n$-form $\phi$ is independent of a K\"ahler metric.
By taking a local section $\psi$ of $K_{\oM}$ non-zero in a neighborhood of $U_i$ of $D$ and
considering $\frac{\phi}{\psi},$ the extension of $\phi$ is reduced to standard result on extension of
$L^2$ holomorphic functions.  The Lemma follows. 

\qed

\ms
We now continue on the proof of Theorem 5.
From Theorem 3, we have 
 the convergence of 
$B^{n,0}_{M_i}(x)$ to $B^{n,0}_{\tM}(x)$ on compacta.  The followings are the two steps that we need to prove
for $i$ sufficiently large,\\
(i)  sections in $H^0_{(2)}(M_i,K_{M_i})$ generate $K_{M_i},$ and\\
(ii) sections in $H^0_{(2)}(M_i,K_{M_i})$ give an immersion of $M_i$.

\ms
Let us first consider (i).  Let $A_i$ be the base locus of
$\Gamma(M_i, K_{M_i}).$  $A_i$ is an algebraic subvariety, which extends to an algebraic variety on a compactification of $M_i$ according to Lemma 1.  Again for simplicity, we denote $K_{M_i}$ by $K.$
$A_i$ is also the set on which the Bergman kernel $B^{n,0}_{M_i}=B^{n,0}_{M_i, K}$ vanishes.  Since $B^{n,0}_{M_i}$
is invariant under biholomorphism, and $M_i\rightarrow M$ is a normal covering, it follows
that $A_i$ is invariant under the deck transformation.  Hence $A_i$ descends to a subvariety
denoted by the same symbol on $M,$ or equivalently, on a fixed fundamental domain $D\subset \tM.$  Moreover, as the covering map is finite, the image on $M$ extends to a subvariety in the compactification $\oM$ of $M$
as well.
Since a $L^2$-holomorphic section of the canonical line bundle on $M_{i}$ pulls back to give
a $L^2$-holomorphic section on $M_{i+1},$ clearly $A_{i+1}\subset A_{i}$ on $M.$  From Noetherian
property, there exists $i_o$ such that $A_i=A_{i+1}$ for $i>i_o.$  Let us denote this set by $A$ 
on $M.$  We are done if $A\cap M=\emptyset,$ which means that the base locus of $K$ on $M_i$ 
is trivial for $i>i_o.$    On the other hand, suppose that $A\cap M\neq\emptyset.$  Let $x\in A\cap M$
and consider a relatively compact neighborhood $U$ of $x$ in $M.$
Pulling it back to the universal covering and denoting by the same symbols, we conclude $B^{n,0}_{M_i}(x)=0$ for all $i>i_o.$  It follows from 
Theorem 2 and 3 that $B^{n,0}_{\tM}(x)=0$ since there is uniform convergence of $B^{n,0}_{M_i}$ to $B_M$ 
on the relatively compact set $U.$  On the other hand, for any bounded domain $\tM$, we can 
always find a non-trivial bounded holomorphic function on $\tM$ non-vanishing at any given
point $x\in \tM$, which means that $B^{n,0}_{\tM}$ is always non-vanishing.  The contradiction establishes base point freeness of $K_{M_i}$ for $i>i_o.$

As for (ii), it is equivalent to showing that there exists $i_o>0$ such that for $i>i_o$
and for any given point $x\in M_i,$ there exist two 
holomorphic sections $f_i, g_i\in \Gamma(M_i, K)$ so
that $d(f_i/g_i)$ is non-degenerate at $x.$  Applying arguments similar to (i) above, in view of the formulations in [Y3], we see readily the validity of (ii).

\qed

\ms
\ni {\bf Remark}   In contrast to the compact Hermitian locally symmetric spaces treated in [Y3], [Y4], since we do not
have a lower bound on the injectivity radius, separation of points by the sections in $H^0_{(2)}(M_i,K_{M_i})$
is not clear.  In particular, the argument in [Y4] comparing heat kernels on $\tM$ and $M_i$ is not
readily applicable.  Additional arguments as in (i) and (ii) above are not sufficient to guarantee
seperation of points on $M_i$ in general.  It can however be proved that for $i$ sufficiently large, 
global sections of $K_{M_i}$ separate distinct points $x, y\in M_i$ except possibly for the case that $\pi_i(x)=\pi_i(y),$ where
 $\pi_i: M_i\rightarrow M$ is the covering map.

\bs
\begin{center} 
{\bf  References}
\end{center}

\ms
\ni [A] Atiyah, M. F., Elliptic operators, discrete groups, and von Neumann
algebras, Ast\'erisque 32-33(1976), 43-72.

\ms
\ni [BM] Barbasch, D., Moscovici, H.,  $L^2$-index and the Selberg trace formula. J. Funct. Anal. 53 (1983), no. 2, 151-201.

\ms
\ni [BW] Borel, A., Wallach, N., Continuous cohomology, discrete subgroups, and representations of reductive groups. Second edition. Mathematical Surveys and Monographs, 67. American Mathematical Society, Providence, RI, 2000.

\ms
\ni
[CG] Cheeger, J., Gromov, M., On the characteristic numbers of complete
manifolds of bounded curvature and finite volume, Differential geometry
and complex analysis, 115-154, Springer, Berlin, 1985.
 
\ms
\ni 
[De] DeGeorge, D., On a theorem of Osborne and Warner. Multiplicities in the cuspidal spectrum. J. Funct. Anal. 48 (1982), no. 1, 81-94.
 
\ms
\ni [DW] DeGeorge, D.,  Wallach, N. R.,  Limit formulas for multiplicities in $L^2(G/\Gamma)$, Ann. of Math. (2) 107 (1978), no. 1, 133-150. 
 
\ms
\ni [Do] Donnelly, H., Elliptic operators and covers of Riemannian manifolds, Math. Zeit. 223 (1996), 303--308.

\ms
\ni [DF]  Donnelly, H., Fefferman, C.,  $L^2$-cohomology and index theorem for the Bergman metric. Ann. of Math. 118 (1983), 593-618.

\ms
\ni [G] Grossman, E. K., On the residual finiteness of certain mapping class groups. J. London Math. Soc. (2) 9 (1974/75), 160-164.

\ms
\ni [HM] Harris, J., Morrison, I., Moduli of curves, Graduate Texts in Mathematics, 187. Springer-Verlag, New York, 1998.

\ms
\ni
[K1] Kazhdan D., On arithmetic varieties, Lie Groups and Their Representations (Proc. Summer School, Bolyai 
J\'anos Math. Soc., Budapest, 1971), Halsted, New York, 1975, pp. 151-217.

\ms
\ni [K2] Kazhdan D., On arithmetic varieties II. Israel J. of Math. 44 (1983), 139Ð159.

\ms
\ni
[L] L\"uck, W., Approximating $L^2$-invariants by their finite-dimensional analogue, Geom. Funct. Anal. 4
(1994) 455Ð481.

\ms
\ni
[Mc] McMullen, C. T., The moduli space of Riemann surfaces is K\"ahler hyperbolic, Ann. of Math.
151 (2000), 327-357.

\ms
\ni
[Mi] Milne, J. S., Kazhdan's Theorem on arithmetic varieties, arXiv:math/0106197v2.

\ms
\ni
[Mo1] Mok, N., Topics in complex differential geometry. Recent topics in differential and analytic geometry, 1Ð141, Adv. Stud. Pure Math., 18-I, Academic Press, Boston, MA, 1990. 

\ms
\ni
[Mo2] Mok, N., Projective-algebraicity of minimal compactifications of complex-hyperbolic space forms of finite volume,
preprint, http://www.hku.hk/math/imr/.

\ms
\ni [MY] Mok, N., Yeung, S-K., Geometric realizations of uniformization of conjugates of Hermitian locally symmetric manifolds. Complex analysis and geometry, 253--270, Univ. Ser. Math., Plenum, New York, 1993.

\ms
\ni [Mu] Mumford, D., Hirzebruch's proportionality theorem in the non-compact case. Invent. Math. 42(1977), 239-272. 

\ms
\ni [NR] Nori, M., Raghunathan, M., On conjugation of locally symmetric arithmetic varieties. 
Proc. Indo-French Conference on Geometry (S. Ramanan and A. Beauville, eds.). Hindustan Book Agency, 1993, 111-122. 

\ms 
\ni [O1]  Ohsawa, T., A remark on Kazhdan's theorem on sequences of Bergman metrics. Kyushu J. Math. 63 (2009), no. 1, 133-137.

\ms 
\ni [O2]  Ohsawa, T., A tower of Riemann surfaces whose Bergman kernels jump at the roof,  Publ. Res. Inst. Math. Sci. 46 (2010), 473-478.

\ms
\ni[PY] Prasad, G., Yeung, S-K., Arithmetic fake projective spaces and arithmetic fake Grassmannians, Amer. J. Math.131(2009), 379-407.

\ms
\ni [R] Rhodes, J., Sequences of metrics on compact Riemann surfaces, Duke Math. J., 72(1993), 725-738.

\ms
\ni [SX] Sarnak, P.,  Xue, X., Bounds for multiplicities of automorphic representations. Duke Math. J. 64 (1991), 207-227.

\ms
\ni
[S] Savin, G., Limit multiplicities of cusps forms, Invent. Math. 95(1989), 149-159. 

\ms
\ni
[SZ] Shiffman, B., Zelditch, S., Distribution of zeros of random and quantum chaotic sections of 
positive line bundles, Comm. Math. Phys. 200 (1999) 661-683.

\ms
\ni[Ti] Tits, J., Classification of algebraic semisimple groups, Algebraic Groups and Discontinuous
Subgroups. Proc.\,A.M.S.\,Symp.\,Pure Math.\: 9(1966) pp. 33--62. 

\ms
\ni
[To] To, W.-K., Distribution of zeros of sections of canonical line bundles over towers
of covers, J. London Math. Soc. (2) 63 (2001), 387-399. 

\ms
\ni
[X] Xue, X. On the Betti numbers of a hyperbolic manifold. Geom. Funct. Anal. 2 (1992), 126-136.

\ms
\ni 

\ms
\ni[Y1] Yeung, S.-K.,
Betti numbers on towers of manifolds, Duke Math. J. 73(1994), 201-226.

\ms
\ni
[Y2] Yeung, S.-K., Properties of complete non-compact K\"ahler surfaces of negative Ricci curvature. J. Reine Angew. Math. 464 (1995), 129-141.

\ms
\ni[Y3] Yeung, S.-K., Very ampleness of line bundles and canonical embedding of coverings of manifolds,
Comp. Math. 123 (2000), 209-223

\ms
\ni[Y4] Yeung, S.-K., Effective estimates on the very ampleness of the canonical line bundle of locally Hermitian symmetric spaces, Tran. AMS, 353(2001), 1387-1401.

\ms
\ni
[Y5] Yeung, S.-K., Geometry of domains with uniform squeezing properties, Adv. in Math., 221(2009), 547-569.

 \end{document}